\documentclass[12pt]{amsart}
\usepackage{etex,stmaryrd}

\makeindex

\usepackage{graphicx,tikz,url,hyperref,epigraph}
\usepackage[utf8]{inputenc} \usepackage[T1]{fontenc} \usepackage{lmodern}
\usepackage{pgf}
\usepackage{epsfig,graphicx,color,amsmath,a4wide,amssymb,amsfonts,amsbsy,xy,latexsym,epsf,pstricks,verbatim,times}
\usepackage{MnSymbol}
\usepackage{pgf,color}
\usepackage{changebar}

\definecolor{LightGrey}{rgb}{.85,.85,.85}
\definecolor{DarkGrey}{rgb}{.5,.5,.5}

\definecolor{Blue}{rgb}{.0,.0,0.9}
\definecolor{LightBlue1}{rgb}{.2,.4,0.9}
\definecolor{LightBlue2}{rgb}{.3,.5,0.9}
\definecolor{LightBlue3}{rgb}{.4,.6,0.9}
\definecolor{LightBlue4}{rgb}{.5,.7,.9}
\definecolor{LightBlue5}{rgb}{.6,.8,.9}
\definecolor{LightBlue6}{rgb}{.7,.9,.9}

\definecolor{Red}{rgb}{.9,.0,.0}
\definecolor{LightRed1}{rgb}{0.9,.2,.4}
\definecolor{LightRed2}{rgb}{0.9,.3,.5}
\definecolor{LightRed3}{rgb}{0.9,.4,.6}
\definecolor{LightRed4}{rgb}{.9,.5,.7}
\definecolor{LightRed5}{rgb}{.9,.6,.8}
\definecolor{LightRed6}{rgb}{.9,.7,.9}

\def\pp{{\mathfrak q}}

\def\KK{{\mathbf K}}
\def\Kb{{ {\bar K}}}

\definecolor{Grey}{rgb}{.5,.5,.5}

\definecolor{Blue}{rgb}{.0,.0,0.9}
\definecolor{LightBlue1}{rgb}{.2,.4,0.9}
\definecolor{LightBlue2}{rgb}{.3,.5,0.9}
\definecolor{LightBlue3}{rgb}{.4,.6,0.9}
\definecolor{LightBlue4}{rgb}{.5,.7,.9}
\definecolor{LightBlue5}{rgb}{.6,.8,.9}
\definecolor{LightBlue6}{rgb}{.7,.9,.9}

\definecolor{Red}{rgb}{.9,.0,.0}
\definecolor{LightRed1}{rgb}{0.9,.2,.4}
\definecolor{LightRed2}{rgb}{0.9,.3,.5}
\definecolor{LightRed3}{rgb}{0.9,.4,.6}
\definecolor{LightRed4}{rgb}{.9,.5,.7}
\definecolor{LightRed5}{rgb}{.9,.6,.8}
\definecolor{LightRed6}{rgb}{.9,.7,.9}

\xyoption{all}
\usepackage[english]{babel}

\newcounter{noalgo}[section]
\newdimen\indentalgo
\newdimen\indentalgodec\indentalgo=0.0mm\indentalgodec=10mm

\newcommand{\If}{\advance\indentalgo by \indentalgodec {\bf if }}

\newcommand{\For}{\global\advance\indentalgo by \indentalgodec {\bf for }}

\newcommand{\Endindent}{\global\advance\indentalgo by -\indentalgodec}

\newdimen\decalage \decalage=0.5cm
\newcounter{algo} 

\newcommand{\PP}{\mathbf P}

\def\<<{\leavevmode
  \raise0.28ex\hbox{$\scriptscriptstyle\langle\!\langle$}\nobreak
  \hskip -.6pt plus.3pt minus.2pt\,}
\def\>>{\,\nobreak\hskip -.6pt plus.3pt minus.2pt
  \raise0.28ex\hbox{$\scriptscriptstyle\rangle\!\rangle$}}

\def\AA{{\mathbf A}}

\def\vol{{{v}}}

\def\Proj{\mathop{\rm{Proj}}\nolimits }
\def\Hom{\mathop{\rm{Hom}}\nolimits }

\def\Sym{\mathop{\rm{Sym}}\nolimits }
\def\IM{\mathop{\rm{Im}}\nolimits }
\def\Ker{\mathop{\rm{Ker}}\nolimits }

\def\Spec{\mathop{\rm{Spec}}\nolimits }

\def\CC{{\mathbf C}}
\def\OO{{\mathbf {O}}}

\def\FF{{\mathbf F}}

\def\AA{{\mathbf A}}

\def\QQ{{\mathbf Q}}

\def\RR{{\mathbf R}}
\def\KK{{\mathbf K}}

\def\ZZ{{\mathbf Z}}
\def\Dgot{{\mathfrak D}}

\def\cC{{\mathcal C}}

\def\cI{{\mathcal I}}

\def\cL{{\mathcal L}}

\def\cM{{\mathcal M}}

\def\cO{{\mathcal O}}

\def\cQ{{\mathcal Q}}

\newtheorem{lemma}{Lemma}

\newtheorem{proposition}{Proposition}
\newtheorem{theorem}{Theorem}

\providecommand{\myproofname}{Proof}

\usepackage{enumitem}

\begin{document}

\begin{abstract}
  We construct   compact descriptions of function  fields
and number fields.
\end{abstract}

\title{Short models of global  fields}

\author{Jean-Marc Couveignes}
\address{Jean-Marc Couveignes, Univ. Bordeaux, CNRS,
  Bordeaux-INP, IMB, UMR 5251, F-33400 Talence, France.}
\address{Jean-Marc Couveignes, INRIA, F-33400 Talence, France.}
\email{Jean-Marc.Couveignes@u-bordeaux.fr}

\date{\today}

\maketitle
\setcounter{tocdepth}{2} 
\tableofcontents

Given a field of scalars $K$ and an indeterminate $x$,
a  natural way to describe a function  field  $K(C)/K(x)$ 
over $K(x)$ with
degree $n$ and genus $g$ is to choose a primitive element  in it 
and to give its minimal polynomial. This is a polynomial
in $K(x)[y]$ and even in $K[x,y]$ once we have cleared
denominators. It has  degree $n$ in $y$ and its degree
in $x$ cannot be expected to be significantly smaller
than the genus $g$ of the curve $C$ in general. So  this plane model involves about $ng$
coefficients while the dimension of the  underlying
modular
problem is rather linear in $n+g$.
Indeed a plane model of a general curve of genus $g$
tends to have   singularities and  tells
more about these singularities than about
the function field itself.
Another more natural possibility to represent   a function field $K(C)/K(x)$
is to 
consider the canonical model of the curve $C$ and express $x$
as a rational fraction of the holomorphic differentials. Even Petri's description of the ideal of 
canonical curves \cite[Chapter III, \S 3]{Arb}  however
involves a number of parameters
that is quadratic in the genus $g$.

Similarly, a natural way to describe a number field
$\KK /\QQ$ of degree $n$ and root discriminant
$\delta_\KK$
is to choose a primitive element  in it 
and give its minimal polynomial. This is a polynomial
in $\QQ [x]$ and even in $\ZZ [x]$ once we have cleared
denominators. It has  degree $n$ in $x$ and its height
cannot be expected to be significantly smaller
than the absolute value of the discriminant $|d_\KK|=\delta_\KK^n$ in general. So  this  model has bit size
of order $n\times \log |d_\KK|$
while the logarithm of the number of number
fields of degree $n$ and
discriminant bounded by $D$ in absolute value
is conjectured \cite{CDO, Malle},
and even proven for $n\leqslant 5$   \cite{DH,Bh4,Bh5}, to be bounded by
$\log D$ plus  a function
of $n$ alone.
The reason for this discrepancy
is that  the quotient  ring of $\ZZ[x]$ by
a unitary irreducible polynomial
tends to be highly singular and provides  more information
than the mere isomorphism class of the corresponding number
field.

In this text we construct 
two families
of natural models for a function field of genus $g$ and degree $n$ over $K(x)$ where $K$ is a field   of cardinality $q$.
These
models involve a number of 
 coefficients in $K$  of  order $n(\log  n)^3(g/n+1+\log_q n)$
 and  $n(\log  n)^2(g/n+1+\log_qn)$ respectively.
 In case $q=+\infty$ we set $\log_q n=0$ in these formulae.
 The first model is obtained by picking few functions
 $(\kappa_j)_{1\leqslant j\leqslant r}$
 of small degree
 in $K(C)$ and 
 finding enough relations with small degree  between them
to produce local equations 
for an affine model
of $C$.
The second model is obtained from the first one
by providing long
enough formal  series  expansions of the $\kappa_j$ to be able
to recover the equations. This second model
is  unambiguous: we can recover the  full ideal
of relations from these formal expansions, not just local
equations. It is also a
bit smaller. It consists of a finite $0$-dimentional affine $K$-scheme (a finite
set in a finite $K$-algebra) that pins down a model of $C$.

Using (not so) parallel constructions in number theory
we propose two families of models
for a number field, having  bit size
of order $n(\log n)^3(\log \delta_K+\log n)$
and  $n(\log n)^2(\log \delta_K+\log n)$
respectively.
The first  model, which was already presented in \cite{jmc2019},
is obtained
by  picking  few  small integers $(\kappa_j)_{1\leqslant j\leqslant r}$
in $\KK$ and finding
enough relations with small degree and small height between them
to produce local equations 
for an affine model
of $\Spec \KK$.
The second model is obtained from the first one
by picking a small unramified prime $p$ an providing long
enough $p$-adic approximations  of the $\kappa_j$ to characterize
the number field. Again this second model
is  unambiguous and a bit
 smaller. It consists of a finite $0$-dimentional affine scheme (a finite
set in a finite algebra) that pins down a model of $\Spec \KK$.

In this text, the
notation $\cQ$ stands for a positive absolute constant.
Any sentence  containing
this symbol becomes true if the symbol is replaced in every occurrence by some 
large enough real number.

\section{Definitions and main statements}\label{sec:main}

We recall notation about function  fields and
number  fields.
We then
define two ways of specifying
a function   field $K(C)$.
We can   provide  local equations for 
a $0$-dimentional variety over the rational field $K(x)$   such
that the function  field $K(C)$
is the residue field of this variety.
We can also  provide formal  expansions  for the
coordinates of some formal point on the former model.
We  state the existence of small models of either kind for any function field $K(C)/K(x)$.
We similarly  define two ways of specifying  
a number  field $\KK$. One option is
to  provide  local equations for 
a $0$-dimentional variety over $\QQ$   such
that the number field $\KK$
is the residue field of this variety.
Another option is to provide $p$-adic expansions  for the
coordinates of some $p$-infinitesimal  point on the former model, with enough accuracy to  recover the equations.
We  state the existence of small models of either kind for any number field $\KK$.

\subsection{Function fields}\label{sec:defff}

Let $K$ be a field.
Let $q$ be the cardinality of $K$.
If $q$ is finite we write $\log_q$ for the logarithm function
with base $q$. If $q=+\infty$ then $\log_q$ is  the constant
zero function.
Let $\PP^1_K =\Proj K[x_0,x_1]$
be the projective line over $K$.
Let  $x=x_1/x_0$ 
and let $\AA_K^1 = \Spec K[x]$
be the affine  line over $K$.
We call $\infty$ the point $(0,1)$ on $\PP^1_K$  with
projective  coordinate $x_0=0$.
Let $C$ be a smooth  absolutely integral projective curve over $K$.
  Let $f : C \rightarrow \PP^1$
  be a finite separable map and let $K(C)/K(x)$ be the corresponding
  regular field extension. 
  In this paper
  we shall call such an extension a {\bf function field}. 
We shall denote   $n$   the degree of  $f$ and 
$g$ the genus of $C$.

\subsection{Function fields as incomplete intersections}\label{sec:ffincom}
Let $K(C)/K(x)$ be a function field as in Section \ref{sec:defff}.
  Let $r\geqslant 1$, $d_x\geqslant 1$  and $d_y\geqslant 1$ be integers.
  Let $y_1$, \ldots, $y_r$ be indeterminates.
  Let  $E_1$, $E_2$, \ldots, $E_r$ be polynomials 
in $K[x][y_1, \ldots , y_r]$ with total degree $\leqslant d_y$ in the $y_j$ and with
degree $\leqslant d_x$ in $x$.
Let $\cI \rightarrow \AA^1=\Spec K[x]$ be the  $K[x]$-scheme with equations
\[E_1 = E_2 = \dots = E_r=0 \text{ and } \det \left( \partial E_i/\partial y_j \right)_{1\leqslant i, \, j\leqslant r} \not = 0.\]
We assume that  the generic fiber $\cI\otimes K(x)\rightarrow  \Spec K(x)$ has an irreducible component
isomorphic to   $\Spec K(C)\rightarrow \Spec K(x)$. We call $\cC
\rightarrow \AA^1 = \Spec K[x]$ the Zariski closure
of this component in $\cI$. We say that $\cC\rightarrow \AA^1$ is an
{\bf incomplete intersection
model} of the function field in dimension $r$ and degrees $(d_x,d_y)$.

\subsection{Interpolant of an infinitesimal   point}\label{sec:interpol}
  Let $r\geqslant 1$, $m\geqslant 1$, $d_x\geqslant 1$ and $d_y\geqslant 1$
  be integers.
  Let $K$ be a field. Let
  $L\supset K$ be a finite separable field extension
  and let $l$ be its degree. Let $\lambda$ be a primitive
  element in $L$.
  Let $S\supset L$
  be a finite separable field extension
  and let $s$ be its degree.
  Let $t$ be an indeterminate.
  Let $b_1^m$, \ldots, $b_r^m$ be $r$ elements in
  $S[[t]]/t^m$.
  Let \[a^\infty=\lambda+t  \in S[[t]] \text{ \,\,  and \,\,}
  a^m=a^\infty  \bmod t^m \in S[[t]]/t^m.\]
  Let
 $x$, $y_1$, \ldots, $y_r$ be $r+1$ indeterminates.
  Let \[\beta^m :  K[x][y_1, \ldots, y_r]
    \rightarrow S[[t]]/t^m\] be 
    the  morphism of $K$-algebras that maps $x$ onto
    $a^m$ and $y_j$ onto $b_j^m$ for $1\leqslant j\leqslant r$.
    We assume that $\beta^m$ is surjective. We denote by 
 $b^m$
    the corresponding $S[[t]]/t^m$ point on 
 $\AA^{r+1}$.

        Let $J_{d_x,d_y}$ be the ideal of $K[x][y_1, \ldots, y_r]$
    generated
  by all polynomials in $\Ker \beta^m$
  having total degree $\leqslant d_y$ in $y_1$, \ldots, $y_r$
  and
  degree $\leqslant d_x$ in $x$.
  We say that a point \[b^\infty =
(a^\infty,b_1^\infty, \ldots, b_r^\infty) \in
(S[[t]])^{r+1}\]
 is the {\bf interpolant} of $b^m$ in degrees
 $(d_x,d_y)$ if
\begin{enumerate}
\item $b^\infty_j \bmod t^m =b^m_j$ for every $1\leqslant j\leqslant r$,
\item all polynomials in  $J_{d_x,d_y}$
 vanish at $b^\infty$,
\item and  $b^\infty$
 is {\bf unique} with the two  above properties.
\end{enumerate}

\subsection{Short models of function fields}\label{sec:shortff}
We shall prove that both
incomplete intersections and interpolants
provide short descriptions of function fields.

\begin{theorem}[Function fields as  small incomplete
    intersections]\label{th:siiff}
  The exists an absolute constant $\cQ$ such that the following is true.
  Let $K$ be a field and let $K(C)/K(x)$ be a
  function field of genus $g\geqslant 2$ and degree $n$
  as in Section \ref{sec:defff}.  Let $q$ be the, possibly
  infinite, cardinality of $K$.
Call $r$ the smallest positive integer such that \[{2r\choose r}\geqslant n(r+1).\] Let $h$ be smallest non-negative integer such that 
  \[q^{h+1}>nr(r+1).\]
    Let
 \[\nu =  2+\lceil \frac{2(g-1)}{n}\rceil  \text{\,\, and  \,\, }d_x=r(h+\nu ) \text{\,\,  and \,\,} d_y=r.\]
  There exists an incomplete intersection model
  of $K(C)/K(x)$ in dimension $r$ and degrees $(d_x,d_y)$.
  The total number of $K$-coefficients in this model
  is \[\leqslant \cQ  (\log  n)^3(g+n(1+\log_qn)).\]
\end{theorem}

The meaning of this theorem is that we have a
short description of $K(C)/K(x)$
as a  quotient  of a finite  algebra : the smooth zero-dimensional
part of a complete intersection of small degree in a projective space of small dimension.
In particular the number of affine parameters required to decribe the incomplete
intersection is almost linear in $g +n$.
There remains a little uncertainty concerning which irreducible
component of the complete  intersection is of interest to us. One way of removing this uncertainty if to specify 
a   geometric point on the targeted component.
We then realize that giving enough such geometric points (counting
multiplicities)
 enables us to reconstruct the equations. This leads us to the  next theorem.

\begin{theorem}[Function fields from interpolants]\label{th:sintff}
  The exists an absolute constant $\cQ$ such that the following is true.
  Let $K$ be a field and let $K(C)/K(x)$ be a
  function field of genus $g\geqslant 2$ and degree $n$
  as in Section \ref{sec:defff}.
Let $q$ be the, possibly
  infinite, cardinality of $K$.
  Let $r$, $h$, $\nu$, $d_x$ and $d_y$ be as in Theorem \ref{th:siiff}.
  Let \[\rho = d_x+r(\nu+h).\]
  There exist  integers $l$, $s$, $m$ and 
 finite separable field extensions $L/K$ and $S/L$ of respective
 degrees $l$ and $s$, and an  element $\lambda$
 in $K$ such that $L=K(\lambda)$, and 
 a point 
\[b^m=(\lambda+t \bmod t^m, b_1^m, \ldots, b_r^m)\in
(S[[t]]/t^m)^{r+1}\] such that
$b^m$ has an
interpolant 
\[b^\infty=(\lambda+t,b_1^\infty, \ldots, b_r^\infty)\in
(S[[t]])^{r+1}\]
in degrees $(d_x,d_y)$. And
there is an  isomorphism of $K$-algebras  from
\[K(\lambda+t,b_1^\infty, \ldots, b_r^\infty)\subset S[[t]]\]
onto 
$K(C)\supset K(x)$ that sends $\lambda+t$ onto $x$. Further 
\[1 \leqslant l \leqslant  \cO(1+\log_q (1+g/n) +\log_qn) \text{\,\, and \,\, } 
1\leqslant s\leqslant n,\] and $m$ is the smallest positive
integer such that $mls>n\rho$, and
\[rmls\leqslant \cQ (\log n)^2\left( n(1+\log_q n)+g\right).\]
\end{theorem}  

The meaning of this theorem is that a model
of the function  field can be recovered from
short formal  expansions of a  few functions in  it.
The  function field is  recovered from
$r$ formal expansions of length $m$ having
coefficients in an extension of degree $ls$
of $K$. This involves $rmls$ coefficients
in $K$. 
So the number of parameters in  this model
is bounded by $\cQ (\log n)^2\left( n(1+\log_q n)+g\right)$.

\subsection{Number  fields}\label{sec:defnf}

By a number field we mean a finite field extension
$\KK/\QQ$ of the rational field. We denote by
$n$ the degree of this extension and by
$\OO$  the ring of integers of $\KK$.
We denote by
$(\rho_i)_{1\leqslant i\leqslant r}$ the $r$ real embeddings of $\KK$
and  by $(\sigma_j, \bar\sigma_j)_{1\leqslant j\leqslant s}$  the $2s$
complex
embeddings of $\KK$.
We also denote by $(\tau_k)_{1\leqslant k\leqslant n}$ the $n=r+2s$
embeddings of $\KK$.
We let \[\KK_\RR = \KK\otimes_\QQ\RR = \RR ^r \times \CC^s\]
be the Minkowski space. We follow the presentation in \cite[Chapitre 1, \S 5]{Neukirch}.
An element $x$ of $\KK_\RR$ can be given  by
$r$ real components $(x_\rho)_\rho$
and $s$ complex components $(x_\sigma)_\sigma$. So we write
$x = ((x_\rho)_\rho , (x_\sigma)_\sigma)$.
For such  an $x$ in $\KK_\RR$ we denote by $||x||$ the
maximum of the absolute values of its $r+s$ components. 
The canonical metric on $\KK_\RR$ is defined
by \[<x,y> = \sum_{1\leqslant i\leqslant r}x_iy_i+
\sum_{1\leqslant j\leqslant s}x_j\bar y_j+\bar x_jy_j.\]
The corresponding Haar measure is said to be canonical also.
The canonical measure  of the convex body $\{x,  ||x||\leqslant 1\}$ is
\[2^r(2\pi )^{s}\geqslant 2^n.\]
The map $a\mapsto a\otimes 1$  injects $\KK$ and $\OO$ into
$\KK_\RR$. For every non-zero $x$ in $\OO$ we have \[||x||\geqslant 1.\]
Let $(\alpha_i)_{1\leqslant i\leqslant n}$ be any $\ZZ$-basis of $\OO$.
Set $A  = (\tau_j (\alpha_i))_{1\leqslant i, j \leqslant n}$. The product $A\bar A^t$
is the Gram matrix $B = (<\alpha_i, \alpha_j>)_{1\leqslant i, j \leqslant n}$ of the canonical form
in the basis $(\alpha_i)_i$.
This is a real symmetric positive matrix.
The volume of $\OO$  according to the canonical Haar measure
is  \[\vol_\OO = \sqrt{\det (B)} = |\det (A)|.\]
We  let 
\[d_\KK =  \det (AA^t)\]
be the discriminant of $\KK$ 
and we denote by \[\delta_\KK= |d_\KK|^{\frac{1}{n}}\]
the root discriminant.
The square of the volume of $\OO$  is  $|d_\KK|$.

\subsection{Number fields as incomplete intersections}\label{sec:nfincom}
Let $\KK/\QQ$ be a number field.
Let $r\geqslant 1$, $d\geqslant 1$, and $H\geqslant 1$
be integers.
Let  $E_1$, $E_2$,
  \dots , $E_r$
  be
  $r$ polynomials
  of total degree $\leqslant d$ in $\ZZ[x_1, \ldots, x_r]$ all having  coefficients bounded in absolute
  value by $H$.
  Let $\cI$ be the  $\ZZ$-scheme with equations
  \[E_1 = E_2 = \dots = E_r=0 \text{ and } \det \left( \partial E_i/\partial x_j \right)_{1\leqslant i, \, j\leqslant r} \not = 0.\] We assume
  that $\cI\otimes \QQ \rightarrow   \Spec \QQ$ 
    contains $\Spec \KK$ as one of its irreducible components.
We call $\cC \rightarrow \Spec \ZZ$  the  Zariski closure
of this component in $\cI$. This is the spectrum of some order
in $\KK$. We say that  $\cC \rightarrow \Spec \ZZ$ is an 
{\bf incomplete intersection
model} of $\KK/\QQ$  in dimension $r$,  degree $d$ and height $H$.

\subsection{Interpolant of a $p$-infinitesimal   point}\label{sec:pinterpol}

  Let $p$ be a prime integer. We fix some algebraic closure
  of $\QQ_p$ and call $\QQ_p^{\rm nr}$ the maximal unramified
  extension in it. For every power $q=p^f$ of $p$
  we call $\QQ_q$ the unique subextension  of $\QQ_p^{\rm nr}/\QQ_p$ of degree $f$
  and $\ZZ_q$ its ring of integers.
Let $r\geqslant 1$, $s\geqslant 1$,
$m\geqslant 1$, $d\geqslant 1$ and $H\geqslant 1$
  be integers.
        Let
    $b_1^m$, \ldots,
        $b_r^m$ be $r$ elements in $\ZZ_{p^s}/p^m$.
Let
    $x_1$, \ldots, $x_r$ be $r$ indeterminates.
    Let \[\beta^m :  \ZZ[x_1, \ldots, x_r]
    \rightarrow \ZZ_{p^s}/p^m\] be 
    the ring homomorphism sending $x_j$ onto
    $b_j^m$ for $1\leqslant j\leqslant r$. We assume
    that $\beta^m$ is surjective.
        We call  $b^m=(b_1^m, \ldots, b_r^m)$
    the corresponding point in $\AA^r(\ZZ_{p^s}/p^m)$.
                Let $J_{H,d}$ be the ideal of $\ZZ[x_1, \ldots, x_r]$
    generated
  by all polynomials in $\Ker \beta^m$
  having total degree $\leqslant d$ 
and all their coefficients bounded by $H$ in absolute value.
We say that a point \[b^\infty =
(b_1^\infty, \ldots, b_r^\infty) \in
(\ZZ_{p^s})^r\]
 is the {\bf interpolant} of $b^m$ in degree
 $d$ and height $H$ if 
\begin{enumerate}
\item $b^\infty_j \bmod p^m =b^m_j$ for every $1\leqslant j\leqslant r$,
\item all polynomials in  $J_{H,d}$
 vanish at $b^\infty$,
\item and  $b^\infty$
 is {\bf unique} with the two  above properties.
\end{enumerate}

\subsection{Short models of number fields}\label{sec:shortnf}
We shall prove that both
incomplete intersections and interpolants
provide short descriptions of number  fields.

\begin{theorem}[Number fields have small incomplete intersection
    models]\label{th:smnf}
  There exists a positive constant $\cQ$ such that the following is true.
  Let $\KK$ be a number field of degree $n\geqslant \cQ$ and
  root discriminant
$\delta_\KK$.
  Let $r$ be the smallest positive integer such that \[{2r \choose r} \geqslant n(r+1).\]
Set
\begin{equation}\label{eq:ell}
  \ell = {2r\choose r}-n \text{ and  }   H = \ell^{\ell/2n}
    \times {2r\choose r}^{1/2}\left(n^2r(r+1)\delta_\KK^{2}\right)^{r}.\end{equation}    Then $r\leqslant \cQ \log n$ and  
 \[\log H\leqslant \cQ{\log n} (\log n+\log \delta_\KK)\] and 
  $\KK$ has an incomplete intersection model in dimension
 $r$, degree $d=r$ and height $H$. Further
\[{d+r\choose r} = {2r\choose r}\leqslant \cQ n\log n\] so that
the total number of coefficients in the $r$ equations is
$\leqslant \cQ n (\log n)^2$.
\end{theorem}

The meaning of Theorem \ref{th:smnf} is that we
have a short description of $\KK/\QQ$
as a quotient of a finite algebra
associated with
the smooth zero-dimensional
part of a complete intersection of small degree and small height in a projective space of small dimension.
The bit size of the model
is  $\leqslant \cQ n(\log n)^3(\log \delta_K+\log n)$.
Every factor in this estimate is awaited but
the $(\log n)^3$.
Another annoyance is the  remaining
little uncertainty about  which irreducible
component of the complete  intersection is of interest to us. One way of removing this uncertainty if to specify 
a   geometric  point on the targeted component.
We then realize that giving enough such geometric points (counting
multiplicities)
enables us to reconstruct the equations. This leads us to the  next theorem.

\begin{theorem}[Number fields from interpolants]\label{th:sintnf}
  The exists an absolute constant $\cQ$ such that the following is true.
  Let $\KK$ be a number field of degree $n\geqslant \cQ$ and
  root discriminant $\delta_\KK$ over $\QQ$.
  Let $r$ 
  and $H$ be as in Theorem \ref{th:smnf}.
  There exist positive integers  $s\leqslant n$,
  and $m$ and a prime integer \[p\leqslant
\cQ  n(\log n)^2 (\log n+\log \delta_\KK)\]
and a point 
\[b^m=(b_1^m, \ldots, b_r^m) \in (\ZZ_{p^s}/p^m)^r\] such that
$b^m$ has an
interpolant 
\[b^\infty=(b_1^\infty, \ldots, b_r^\infty)\in (\ZZ_{p^s})^r\]
in degree
$r$ and height $H$. And
the field \[\QQ(b_1^\infty, \ldots, b_r^\infty)\subset \QQ_q\]
is isomorphic to $\KK$. And
\[rms\log p\leqslant \cQ n(\log n)^2\left(\log n+\log \delta_\KK \right).\]
\end{theorem}  

The meaning of this theorem is that a model
of the number field can be recovered from
a few short $p$-adic expansions of a few algebraic numbers  in it.
The  number field is  recovered from
$s$, $m$, $p$ and by
$r$ elements in $\ZZ_q/p^m\ZZ_q$ where $q=p^s$.
Each of these $r\leqslant \cQ\log n$ elements is a  $q$-adic expansion of bit size
$ms\log p$.
So the bit size of this model
is bounded by
$\cQ  n(\log n)^2\left(\log n+\log \delta_\KK \right)$.

\section{Constructing  models of function fields}\label{sec:fff}

In this section we prove Theorems \ref{th:siiff} and \ref{th:sintff}.
In Section \ref{sec:shortfu} we recall the definition
of the Maroni invariants and we bound them.
In Section \ref{sec:smalleqff} we use a  generic well-posedness
theorem in multivariate Hermite interpolation 
due to  Alexander
and   Hirschowitz in order to prove the existence
of  small affine models of function fields.
In Section \ref{sec:interpolff} we prove that
formal expansions of coordinates at a smooth point
in these models characterize the function field.

\subsection{Small degree functions}\label{sec:shortfu}

Let $K$ be a field. Let $K(C)/K(x)$
be a function field as in Section \ref{sec:defff}.
We call $C_{\rm aff}$ the affine part
\[C_{\rm aff}=f^{-1}(\AA^1)\subset C\] of $C$.
The sheaf $f_\star\cO_C$ is a rank $n$ locally free  ${\PP^1}$-module.
Every such sheaf  decomposes as a direct sum of invertible
sheaves. So 
\begin{equation}\label{eq:decomp}
f_\star\cO_C = \bigoplus_{0\leqslant i \leqslant  n-1}\cO_{\PP^1}(-a_i)\end{equation}
where $a_0=0$ and the remaining $a_i$ form a non decreasing sequence of
strictly positive integers, called the scrollar invariants,
or the Maroni invariants of the cover. We call $\cO_i$ the $i$-th term on the right
hand side of Equation~(\ref{eq:decomp}).
Let $\tilde \omega_i$ be a generator of $H^0(\PP^1,\cO_i \otimes \cO_{\PP^1}(a_i))$ and set
\[\omega_i={x_0^{-a_i}}{\tilde\omega_i}.\]
This is a function  on $C$, regular on $C_{\rm aff}$,  and  the maximum order of its  poles above  $\infty$  is  $a_i$.
The family $(\omega_i)_{0\leqslant i\leqslant d-1}$
is a $K[x]$ basis of \[H^0(\AA\!\!{}^1,f_\star (\cO_{C})) = K[C_{\rm aff}]\] the integral closure of $K[x]$
in $K(C)$.
The determinant sheaf \[\det(f_\star \cO_C)=\wedge^df_\star \cO_C\] is isomorphic to
$\cO_{\PP^1}(-a)$ where \[a = \sum_{0\leqslant i\leqslant n-1}a_i = n+g-1.\]
Because $f_\star \cO_C$
is a sheaf of ${\PP^1}$-algebras, there exists
a multiplication tensor \[ m \in H^0(\PP^1,\widehat{f_\star \cO_C}
\otimes \widehat{f_\star \cO_C} \otimes f_\star \cO_C).\]
For every triple
$(i,j,k)$ of integers with 
$0\leqslant i,j,k \leqslant n-1$, we call
\[m_{i,j}^k \in H^0(\PP^1,\hat \cO_i\otimes \hat
\cO_j\otimes \cO_k)\]
the corresponding component of $m$. Since
$m_{i,j}^k$ is a global section of a sheaf isomorphic to
$\cO_{\PP^1}(a_i+a_j-a_k)$ we have 
\begin{equation}\label{eq:ineqa}
  a_i+a_j\geqslant a_k\end{equation} whenever $m_{i,j}^k$ is non-zero.
We call $\mu_{i,j}^k\in K[x]$ the coefficients of the
multiplication table of $K[C_{\rm aff}]$ in the basis
$(\omega_0, \ldots, \omega_{n-1})$. We
have 
\[\omega_i\omega_j=\sum_{k}\mu_{i,j}^k\omega_k\]
and $\mu_{i,j}^k$ is a polynomial of degree $\leqslant a_i+a_j-a_k$.
In particular $\mu_{i,j}^k$ is zero unless
inequality (\ref{eq:ineqa}) holds true.

There exists
a permutation $\sigma$ of the set of
integers in the interval $[1,n-2]$ such that,
for every integer  $i$ in this interval,
the coefficient $\mu_{i,\sigma(i)}^{n-1}$ is non zero.
Otherwise the determinant $\det (\mu_{i,j}^{n-1})_{1\leqslant i,j\leqslant n-2}$
would be zero and there would exist a function $\eta$ in the $K(x)$-vector space  generated
by  $\omega_1$, $\omega_2$, \ldots, $\omega_{n-2}$ such that $\eta W\subset W$ where $W$ is the
$K(x)$-vector space  generated
by  $\omega_0=1$, $\omega_1$, \ldots, $\omega_{n-2}$. This would turn $W$ into a $K(x,\eta)$-vector space. But this is imposible
because the extension $K(x,\eta)/K(x)$ has degree 
at least two and the codimension of $W$ in $K(C)$ is $1$.
We deduce that $a_i+a_{\sigma (i)}\geqslant a_{n-1}$ for $1\leqslant i\leqslant n-2$
and summing out  we find \[ (n-2)a_{n-1} \leqslant 2\sum_{1\leqslant i\leqslant n-2}a_i  = 2(a-a_{n-1})\]
that is \begin{equation}\label{eq:maroni}
  a_{n-1}\leqslant \frac{2a}{n} = 2\left( 1+\frac{g-1}{n}\right).\end{equation}
This inequality and its proof are  the function field analogue  to
\cite{Bha}[Theorem 3.1].
\begin{lemma}
The Maroni invariants of a function field
of degree $n$ and genus $g$ are bounded form above
by
\[\nu  
= 2+\Bigl\lceil \frac{2(g-1)}{n}\Bigr\rceil\]
the smallest integer bigger than or equal to $2+2(g-1)/n$.
\end{lemma}
So the field
extension $K(C)/K(x)$ has a basis made of functions
in $K[C_{\rm aff}]$ with  small degree.

\subsection{Small equations for function fields}\label{sec:smalleqff}

Having found small degree functions $(\omega_i)_{0\leqslant i\leqslant n-1}$
we now look for small degree
relations between them.
Let $r\geqslant 2$ and $h\geqslant 0$ be integers. For every
pair of integers  $(i,j)$ with $0\leqslant i\leqslant n-1$ and
$1\leqslant j\leqslant r$ we let \[u_{i,j}(x)\in K[x]_h\]
be a polynomial of degree $\leqslant h$ in $x$.
For $1\leqslant j\leqslant r$ we set
\[\kappa_j = \sum_{0\leqslant i\leqslant n-1}u_{i,j}\omega_i \in K[C_{\rm aff}].\]
The functions $\kappa_j$ have all their poles above $\infty$ and the order
of these poles is bounded from above by $h+\nu$.
We define a morphism of $K[x]$-algebras 
\[\epsilon_{\rm aff} :  K[x][y_1, \ldots, y_r]\rightarrow K[C_{\rm aff}]\] by sending $y_j$ to $\kappa_j$ for $1\leqslant j\leqslant r$.
This  results in  a  morphism 
\[
\xymatrix{e_{\rm aff} \,\,\,\,\,\,\,\,\, : &  C_{\rm aff} \ar@{->}[rrrr]  \ar@{->}_(.4){f_{\rm aff}}[drr] &&&&\AA^{r+1}_K
  =\Spec K[x][y_1, \ldots, y_r]
  \ar@{->}[dll]\\&
  &&\AA\!{}^1_K =\Spec K[x]&&
}
\]
The image by $e_{\rm aff}$ of the cycle $[C_{\rm aff}\otimes_{f_{\rm aff}} K(x)]$ is a $0$-cycle $P$
of degree $n$ on  $\AA^r_{K(x)}$. We let $I$ be the corresponding ideal
of $K(x)[y_1, \ldots, y_r]$ and denote by $2P\subset  \AA^r_{K(x)}$ the subscheme associated with
$I^2$.
We denote by $\cC\rightarrow \AA^1$ the Zariski closure of $P/K(x)$
in $\AA^{r+1}$.
Let $d\geqslant 1$ be an integer such that
\begin{equation}\label{eq:bino}n(r+1)\leqslant
{d+r\choose d}.\end{equation} Let
\[V =K(x)[y_1, \ldots, y_r]_d\] be the $K(x)$-vector
space of polynomials having total degree $\leqslant d$ in the
variables $y_1$,
\ldots, $y_r$ and let 
$M$ be the basis  of $V$ consisting of monomials.
We say that $2P$ is {\bf well poised} in degree $d$
if the restriction to
$V$ of the quotient map \[K(x)[y_1, \ldots, y_r]\rightarrow
K(x)[y_1, \ldots, y_r]/I^2\] is surjective.
We let $\Omega$ be a separable  closure of $K(x)$ and we call
$T$ the set of $K(x)$-embeddings \[\tau : K(C)\rightarrow
\Omega .\]  There are $n$ such embeddings. For every $\tau$ in
$T$ we call \[P_\tau  = (\tau (\kappa_j))_{1\leqslant j\leqslant r}\in \AA^r(\Omega)\] the corresponding  geometric point
of $P$. 
The scheme $2P$ is well poised if and only if
the matrix
\begin{equation*}
\cM_P^1 = [(m(P_\tau))_{ \tau \in T,  \, m\in M  },
  (\partial m/\partial y_1 (P_\tau))_{ \tau \in T,\, m\in M},
  \dots, 
  (\partial m/\partial y_r (P_\tau))_{\tau \in T, \, m\in M}]
\end{equation*}
with $n(r+1)$ rows  and ${d+r\choose d}$ columns
has maximal rank $n (r+1)$.
We note
that $\cM_P^1$ consists of $r+1$ blocks of size
$n\times {d+r\choose d}$ piled vertically.
It
has maximal rank  for a generic $P$ when
$d\geqslant 5$, according to the theorem
of Alexander  and Hirschowitz \cite{Alexander, AH}. 

    The maximal minors of  $\cM_P^1$ are polynomials of  degree
    $\leqslant dn(r+1)$ in each of
    the $u_{i,j}$ and one of them is not identically zero.
    The latter determinant  cannot vanish on the cartesian product
    $(K[x]_h)^{nr}$ as soon as the cardinality
    of $K[x]_h$ is bigger than $dn(r+1)$. 
If $K$ is a finite field
    with cardinality $q$ this condition is granted as soon as \[q^{h+1}>
    dn(r+1).\]If $K$ is  infinite we can
afford $h=0$. 
We will assume  that  $2P$ is well
poised in degree $d$. So $P$ is well poised also and the map $e_{\rm aff}\otimes K(x)$
is a closed immersion.

We look for small degree equations between the $(\kappa_j)_{1\leqslant j\leqslant r}$. More precisely we look for 
polynomials in \[\Ker \epsilon_{\rm aff} \subset K[x][y_1, \ldots, y_r]\] having total degree $\leqslant d$
in the $y_j$ and  degree in $x$ as small as possible.
We denote $K[x][y_1, \ldots, y_r]_d$ the $K[x]$-module of polynomials with total degree $\leqslant d$
in the $y_j$. This is a free $K[x]$-module of rank ${d+r\choose d}$.
We define a morphism of $K[x]$-modules 
\[K[x][y_1, \ldots, y_r]\rightarrow K[C_{\rm aff}]\otimes x_0^{d(h+\nu)}= H^0(\AA^1,f_\star \cO_C\otimes_{\PP^1} \cO_{\PP^1}(d(h+\nu)))\] by
sending $\prod_jy_j^{\gamma_j}$ to $x_0^{d(h+\nu)}\prod_j\kappa_j^{\gamma_j}$ 
for all exponents  $(\gamma_1, \ldots, \gamma_r)$ with $\sum_{j}\gamma_j\leqslant d$.
This morphism of $K[x]$-modules  extends to a morphism of  $\PP^1$-modules
\[\epsilon_{\leqslant d} : \Sym^{\leqslant d}  \cO_{\PP^1}^{\oplus r} \rightarrow  f_\star \cO_C \otimes_{\PP^1} \cO_{\PP^1}(d(h+\nu))\]that is generically surjective. 
The kernel of $\epsilon_{\leqslant d}$ is a locally free  sheaf
of rank
\[\ell = {d+r\choose d} -n.\] The image of  $\epsilon_{\leqslant d}$ is a locally free  sheaf
of rank $n$ and its determinant has degree \[\deg \det \IM \epsilon_{\leqslant d} \leqslant
\deg \det f_\star\cO_C \otimes_{\PP^1}\cO_{\PP^1}(d(h+\nu))=
-a+nd(h+\nu).\]
Since \[\det \Sym^{\leqslant d} \cO_{\PP^1}^{\oplus r} \simeq
\det \Ker \epsilon_{\leqslant d} \otimes_{\PP^{1}}\det \IM \epsilon_{\leqslant d}\]
we deduce that
\[\deg \det \Ker \epsilon_{\leqslant d} \geqslant
-n(h+\nu)d+a.\]
There exists $\ell$  non negative integers $(e_i)_{1\leqslant i\leqslant \ell}$
that form a non decreasing sequence and such that 
\[\Ker \epsilon_{\leqslant a} \simeq \oplus_{1\leqslant i\leqslant \ell}\cO_{\PP^1}(-e_i).\]
More explicitly there exist  $\ell$ polynomials
$(E_i)_{1\leqslant i\leqslant \ell}$ in $K[x][y_1, \ldots, y_r]_d$
such that the coefficients
of $E_i$ are polynomials in $K[x]$ of degree $\leqslant e_i$ and
\[E_i(\kappa_1, \ldots, \kappa_r)=0.\] Further 
 \[\sum_{1\leqslant i\leqslant \ell}e_i \leqslant    n(h+\nu)d-a.\]

 On the one hand \[ e_i\leqslant \lfloor \frac{n(h+\nu)d-a}{n} \rfloor
 \leqslant (h+\nu)d\]
for every $1\leqslant i\leqslant \ell +1 -n$.
On the other hand the scheme $2P$ is well poised and the  $\Omega$-vector space generated
by the $E_i$ for
$1\leqslant i\leqslant \ell+1-n$ has codimension $n-1 < n$ in $\Ker \epsilon_{\leqslant d} \otimes_{\PP^1}\Omega$. 
So there exists at least one embedding $\tau \in T$ such that the $(\ell +1 -n)\times r$ matrix
\[\left((\partial E_i /  \partial y_j)(P_\tau)\right)_{1\leqslant i\leqslant \ell+1-n, \,  1\leqslant j\leqslant r}\]
has maximal rank $r$. 
This means that there exist $r$ integers
$1\leqslant i_1 < i_2 < \dots < i_r \leqslant \ell +1-n$
such that the minor determinant 
\begin{equation}\label{eq:Phiff}
  \Phi = \det \left(  (\partial E_{i_k} /
  \partial y_j)\right)_{1\leqslant k,\,  j \leqslant r}\end{equation}
does not vanish
at $P_\tau$
for some $\tau$ and thus for all $\tau$ by Galois action.

We now choose $d$, $r$, and  $h$ depending on $g$, $n$, and $q$.
We take $r=d$ and in order to grant the condition
in Equation (\ref{eq:bino}) we choose $r$ to be the smallest integer
such that 
\[\frac{1}{r+1}{2r\choose r}\geqslant n.\]
Since
${2k \choose k}\geqslant 2^k$ for every integer $k\geqslant 1$ we have
\[\frac{1}{k+1}{2k \choose k}\geqslant 2^{\frac{k}{2}}\]
for $k$ large enough.
Further  \[\frac{1}{k+2}{2k+2 \choose k+1}\leqslant \frac{1}{k+1}{2k \choose k}\times 4\]
for $k\geqslant 1$.
We deduce  that \begin{equation}\label{eq:bornea}
  \frac{1}{r+1}{2r\choose r}\leqslant 4n\text{\,\,\,  and  \,\,\, } r\leqslant \cQ \log n.\end{equation}
If $K$ is infinite we set $h=0$. If $K$ has $q$ elements  we take $h$ to be the smallest integer such that
\[q^{h+1}>nr(r+1).\] We check that
\begin{equation}\label{eq:borneh}
  h\leqslant \cQ(1+\log_q n).\end{equation}

We have \[n+\ell = {2r\choose r} \leqslant 4n(r+1)\]
and \begin{equation}\label{eq:bornee}
  e_i\leqslant    d(h+\nu) 
\leqslant \cQ (\log  n)(1+\log_qn+\frac{g}{n})\end{equation}
for every $1\leqslant i\leqslant \ell +1 -n$. The number of
$K$-coefficients in $E_i$ for each $1\leqslant i\leqslant \ell +1 -n$
is thus bounded from above by
\[{2r\choose r}\times (1+\max_{1\leqslant i\leqslant \ell +1 -n }e_i)
\leqslant \cQ  (\log  n)^2(g+n(1+\log_qn)).\]
This finishes the proof of Theorem \ref{th:siiff}.

\subsection{Interpolants for function fields}\label{sec:interpolff}

We
set \[d_x =  (h+\nu)r \]
an upper bound for the degree in $x$
of the  equations $E_{i_k}$.
In order to prove Theorem \ref{th:sintff} we look
for a  fiber of $f_{\rm aff} : C_{\rm aff} \rightarrow \AA^1$
where the determinant $\Phi$ of Equation (\ref{eq:Phiff})
does not vanish. Let $\Psi(x) \in K[x]$
be the product of all $\Phi (P_\tau)$
for $\tau \in T$. We bound the degree of $\Psi(x)$.
First
$\Phi$ has degree $\leqslant rd_x$ in $x$
and  total  degree $\leqslant r(r-1)$ in the $y_j$. Since the $\alpha_j$
have poles of order at most $\nu+h$ we deduce that
the evaluation of $\Phi$ at $(\alpha_j)_{1\leqslant j\leqslant r}$
is a function in $K[C_{\rm aff}]$ with poles of order
$\leqslant (\nu+h)r(r-1)+rd_x$. The norm of this function is
$\Psi(x)$ and its degree in $x$ is
bounded from above by $n$ times $(\nu+h)r(r-1)+rd_x$. So 
\[\deg \Psi \leqslant \cQ n(\log n)^2(1+\log_q n+\frac{g}{n})\] according to
Equations (\ref{eq:bornea}), (\ref{eq:bornee}) and  (\ref{eq:borneh}).

If $K$ is infinite there exist infinitely many unitary
polynomials of degree one in $K[x]$ that are prime to $\Psi (x)$.
If $K$ is finite with cardinality  $q$ then we use the
following lemma.

\begin{lemma}\label{lem:psi}
  Let $K$ be a finite field with $q$ elements.
  Let $\Psi (x)\in K[x]$ be a polynomial with degree
  $\deg \Psi \geqslant 1$.
      There exists
  an irreducible  polynomial in $K[x]$ that is prime to $\Psi$
  and  has degree $\leqslant \log_q (\deg\Psi )+1$.
\end{lemma}

For every integer $k\geqslant 1$ we set
$\Pi_k (x) = x^{q^k}-x$. This is a separable
polynomial in $K[x]$ and all its irreducible factors
have degree dividing $k$. If $q^k> \deg \Psi$
then at least one of these irreducible
factors is prime to $\Psi$.  We take $k$
to be the smallest integer bigger than $\log_q\deg\Psi$.
    \hfill $\Box$

    We apply Lemma \ref{lem:psi} to our $\Psi (x)$ and find
    a unitary irreducible  polynomial  $F(x)\in K[x]$ that is
    prime to $\Psi (x)$ and has degree \[l \leqslant \log_q (\deg\Psi )+1
    \leqslant  \cO(1+\log_q (1+g/n) +\log_qn) .\]
    Let $\Kb$ be an algebraic  closure of $K$. Let $\lambda \in\Kb$
    be a root of $F$. Let $\sigma$ be a point in  $C(\Kb)$
    such that $f(\sigma)=\lambda$. We set  $L = K(\lambda)$
    and $S=K(\sigma)$ the residual field at $\sigma$.
    Let 
    $s$  be the degree
    of $S/L$. We set
    \[\rho = d_x+r(\nu+h) \leqslant \cQ (\log  n)(1+\log_qn+\frac{g}{n})\] and let $m$ be the smallest positive integer
    such that $m l s > n\rho$.
We denote
$b_1^\infty$, \ldots, $b_r^\infty$ the expansions of the functions 
of $\kappa_1$, \ldots, $\kappa_r$ at  the point $\sigma$ in the
local parameter $t=x-\lambda$.
We set \[b_j^m = b_j^\infty \bmod t^m \in  S[[t]]/t^m\]
for every $1\leqslant j\leqslant r$.
The point \[b^m=(\lambda+t \bmod t^m, b_1^m, \ldots, b_r^m)\]
in $\AA^{r+1}(S[[t]]/t^m)$ has an interpolant in degrees $(d_x,d_y)$
and this interpolant is
\[b^\infty = (\lambda +t, b_1^\infty, \ldots, b_r^\infty).\]
Indeed let $E(x,y_1, \ldots, y_r)$ be a polynomial
in $K[x][y_1, \ldots, y_r]$ with degree $\leqslant d_x$
in $x$ and with total degree $\leqslant d_y=r$ in the $y_j$.
The polynomial  $E$
induces a function in $K[C_{\rm aff}]$
having at most
$n$ poles with order $\leqslant d_x+r (\nu +h)=\rho$ each.
If this function vanishes at $b^m$
then the divisor of its zeros
has degree $m l s  > n\rho$.
So $E$ 
induces the zero function on  $C_{\rm aff}$. Therefore
it
vanishes at $b^\infty$.
Unicity follows from the existence of $r$ equations
with degree $\leqslant d_x$ in $x$ and total degree
$\leqslant d_y$ in the $y_j$  having a non-zero Jacobian
determinant at $\sigma$.
Since \[\frac{n\rho}{ls}<m\leqslant \frac{n\rho}{ls}+1\]
we have $mls\leqslant n\rho+ls
\leqslant \cQ (\log  n)(n(1+\log_qn)+{g})$
and $r\leqslant \cQ\log n$.
This finishes the proof of Theorem \ref{th:sintff}.

\section{Constructing models of number fields}\label{sec:nf}

We prove Theorems \ref{th:smnf} and \ref{th:sintnf}.
A natural counterpart to sheaves of modules over $\PP^1$
in the context of number fields are
euclidean lattices.
Sections \ref{sec:smallnf} and
\ref{sec:smalleqnf} summarize \cite{jmc2019}.
These sections  realize  the spectrum of some order
of  $\KK$ as an irreducible
component of a  complete intersection in some affine
space over $\ZZ$.
Section \ref{sec:interpolnf} finds a smooth $\FF_p$-point
on this irreducible component and thickens it just enough
to  characterizes the number field.

\subsection{Small elements in number fields}\label{sec:smallnf}

The ring of integers of a number field can be seen as  a euclidean lattice.
The existence
of an  integral domain structure which is compatible with the $L^\infty$-norm restricts the possibilities
for the successive minima: they must all have the same order of magnitude.
\begin{lemma}[Number fields have small integers]\label{lem:balanced}
  The ring of integers $\OO$ of a number field $\KK$ with degree $n$ and
  root discriminant $\delta_\KK$ contains $n$ linearly independant elements
  $(\omega_i)_{1\leqslant i\leqslant n}$ over $\ZZ$ such that all the absolute
  values of all the $\omega_i$ are $\leqslant \delta_\KK^{\, 2}$.
\end{lemma}

This is Proposition 1 in \cite{jmc2019}.
Bhargava, Shankar,  Taniguchi,  Thorne,  Tsimerman, and  ZhaoSee
prove in \cite{Bha}[Theorem 3.1]  a similar statement.

\subsection{Small equations for number fields}\label{sec:smalleqnf}

Having found small integers
 $(\omega_i)_{0\leqslant i\leqslant n-1}$
we now look for
relations between them with 
small degree and small height.
Let $d\geqslant 5$ and $r\geqslant 1$ be two integers. We assume that
\[n(r+1)\leqslant {d+r\choose d}.\]
Let $M$ be the set of  monomials of total degree $\leqslant d$
in the $r$ variables $x_1$, \ldots, $x_r$.
We have  \[\AA^r_\CC = \Spec \CC[x_1, \ldots, x_r]
\subset \Proj \CC[x_0, x_1, \ldots, x_r] = \PP^r_\CC.\]
Let $V_\CC$ be the $\CC$-linear space generated by $M$.
We may associate to every element in $M$ the corresponding
degree $d$ monomial in the $r+1$ variables
$x_0$, $x_1$, \ldots, $x_r$.
We thus identify $V_\CC$ with $H^0(\cO_{\PP^r_\CC}(d))$, the space  of homogeneous polynomials of degree $d$.
We call $T$ the set of all $n$ embeddings of $\KK$
into $\CC$.
       We let  $(\omega_i)_{1\leqslant i\leqslant n}$ be $n$ linearly independant
    short elements in $\OO$ as in Lemma \ref{lem:balanced}.
We pick $rn$ rational integers $(u_{i,j})_{1\leqslant i\leqslant n, \, 1\leqslant j \leqslant r}$ and we
set
\[\kappa_j =  \sum_{1\leqslant i\leqslant n}u_{i,j}\omega_{i}\]
for $1\leqslant j\leqslant r$. 
For every ${\tau \in T}$ we  consider  the points  \[P_\tau = (        \sum_{1\leqslant i\leqslant n}u_{i,j}\tau(\omega_{i})       )_{1\leqslant j\leqslant r}\in \CC^r = \AA^r(\CC)\subset \PP^r(\CC)\] and call $P$ the union of these $n$ points.
The maximal minors of the matrix
\[\cM_P^1 = [(m(P_\tau))_{ \tau ,  \, m\in M  },
  (\partial m/\partial x_1 (P_\tau))_{ \tau ,\, m\in M}, \dots, 
  (\partial m/\partial x_r (P_\tau))_{\tau , \, m\in M}]
\] 
are polynomials of
degree $\leqslant dn(r+1)$ in each of
the $u_{i,j}$ and one of them is not identically zero
according to the  theorem
    of Alexander and Hirschowitz \cite{Alexander,AH}.     
The latter determinant  cannot vanish on the cartesian product $[0,dn(r+1)]^{nr}$. Thus there exist $nr$ rational  integers $u_{i,j}$ in the range
\[[0,dn(r+1)]\] such that the corresponding scheme $2P$ is well poised:
it imposes $n(r+1)$ independent conditions on degree $d$ homogeneous polynomials. 
We   assume that  the $u_{i,j}$ meet these conditions.
We look for polynomials with degree
 $\leqslant d$ and small integer coefficients vanishing at $P$.
 We denote by $V_\RR = \RR[x_1, \ldots, x_r]_d$  the $\RR$-vector space
 of polynomials in $\RR[x_1, \ldots, x_r]$
 of degree $\leqslant d$. There is a unique $\RR$-bilinear
 form on $V_\RR$ that turns the set $M$
 of monomials into an orthonormal basis.
 Let   $V_\ZZ = \ZZ [x_1, \ldots, x_r]_d$.
 The lattice of relations
 with integer coefficients and degree $\leqslant d$ is  a free $\ZZ$-module $\cL \subset   V_\ZZ$
 of  rank
 \[\ell = {d+r\choose d} -n.\] We set $L = \cL\otimes_\QQ \RR$
 the underlying  $\RR$-vector space
 and $L^\perp$ its orthogonal complement in $V_\RR$. We denote by $\cL^\perp$
 the intersection $\cL^\perp  = L^\perp \cap V_\ZZ$.
 Since $V_\ZZ$ is unimodular,
 $\cL$ and $\cL^\perp$ have the same
 volume. See  \cite[Corollary 1.3.5.]{Martinet}.
 We denote by $\hat \OO = \Hom (\OO,\ZZ)$ the dual
 of $\OO$, the ring of integers of $\KK$,  as a $\ZZ$-module.
The evaluation map
at $(x_1, \ldots, x_r)=(\kappa_1, \ldots, \kappa_r)$ is denoted \[\epsilon_{\ZZ, d} :
 \ZZ [x_1, \ldots, x_r]_d\rightarrow \OO.\]
 We observe that $\cL^\perp$ contains
 the image of $\hat \OO$ by the transpose
 map \[\hat \epsilon_{\ZZ, d} : \hat \OO \rightarrow
 \ZZ [x_1, \ldots, x_r]_d\]
 \noindent where we have identified $\ZZ [x_1, \ldots, x_r]_d$
 with its dual thanks to the canonical bilinear form.
 So the volume of $\cL$ is bounded from above by
 the volume of  $\hat \epsilon_{\ZZ, d} (\hat \OO)$.
 We consider the matrix
 \[\cM_P^0 = [(m(P_\tau))_{\tau , \, m\in M }]\]
 of the map $\epsilon_{\CC, d} = \epsilon_{\ZZ, d}\otimes_\ZZ\CC$
 in the canonical bases.
 If we prefer to use an integral basis of $\OO$ on the right
 we should multiply $\cM_P^0$ on the left 
 by the inverse $T$ of the matrix  of a basis of $\OO$ in the
 canonical basis. We deduce that the square of  the volume
 of $\hat \epsilon_{\ZZ, d} (\hat \OO)$ is the determinant
 of $T\cM_P^0(\cM_P^0)^tT^t$.
 Since $T\cM_P^0$ has real coefficients we have
 \[\det (T\cM_P^0(\cM_P^0)^tT^t) = \det \left(T\cM_P^0\left( \overline{\cM_P^0}\right)^t\bar T^t\right)
 =\det \left(\cM_P^0   \left( \overline{\cM_P^0}\right)^t  \right) / |d_\KK|.\]
 So the square of the volume of the lattice of relations is bounded by
 the determinant of the hermitian positive definite matrix
 $\cM_P^0\left( \overline{\cM_P^0}\right)^t$.
 We will use several times the following lemma.
 \begin{lemma}\label{lemma:eval}
   If $E(x_1, \ldots, x_r)\in \ZZ [x_1, \ldots, x_r]$ is a polynomial
   with degree $\leqslant d_E$ and height $H_E$ then for every embedding
   $\tau$
   the evaluation
   of $E$ at $P_\tau$ is bounded in absolute value
   by \[H_E\left(n^2d(r+1)\delta_\KK^{2}\right)^{d_E}{d_E+r\choose r}.\]
\end{lemma}   

 The coefficients $u_{i,j}$ are bounded 
 by $dn(r+1)$. All the absolute values of the $\omega_i$ are bounded
 by $\delta_\KK^{2}$. So the terms
 in $E$
 are bounded from above by
 \[H_E\left(n^2d(r+1)\delta_\KK^{2}\right)^{d_E}\] and there
 are at most ${d_E+r\choose r}$ of them.
 \hfill $\Box$

 Recall that the coefficients in $\cM_P^0$
 are degree $\leqslant d$ monomials in the $\kappa_j$.
 We deduce from
 Lemma \ref{lemma:eval}
 that they are 
bounded from above by \begin{equation*}
   \left( n^2d(r+1) \delta_\KK^2\right)^d.\end{equation*}
 The coefficients in $\cM_P^0\left( \overline{\cM_P^0}\right)^t$ are bounded from
 above by \[\Dgot = {d+r\choose d}\left( n^2d(r+1)\delta_\KK^{2}\right)^{2d}.\]
 The  matrix $\cM_P^0\left( \overline{\cM_P^0}\right)^t$
 being hermitian positive
 definite, its determinant is bounded
 from above by the product of the diagonal terms.
 We deduce that the volume of the lattice $\cL$
 of relations  is bounded from above by $\Dgot ^{n/2}$.
Recall that the dimension of $\cL$ is \[\ell = {d+r\choose d}-n.\]
For any $x$  in $V_\RR$ we denote by $||x||$ the $\ell_2$-norm in
the monomial basis. The volume of the sphere  $\{x \in L,  ||x||\leqslant 1\}$ is $\geqslant 2^{\ell}\ell^{-\ell/2}$.
Applying Minkowski's second theorem \cite[Lecture III, \S 4, Theorem 16]{Siegel}
to the gauge function $x\mapsto ||x||$ we find that $\cL$ contains $\ell$ linearly independant
elements $E_1$, $E_2$, \ldots, $E_\ell$ such that
\[\prod_{1\leqslant i\leqslant \ell}||E_i||\leqslant \ell^{\ell/2}\Dgot^{n/2}.\]
We assume that the sequence $i\mapsto ||E_i||$
is non-decreasing and deduce that 
\[||E_i||\leqslant \ell^{\frac{\ell}{2(\ell+1-i)}} \Dgot^{\frac{n}{2(\ell+1-i)}}\]for every $1\leqslant i\leqslant \ell$.
We forget
the last $n-1$ equations.
On the one hand \[ ||E_i||\leqslant \ell^{\ell/2n} \Dgot^{1/2}\]
for every $1\leqslant i\leqslant \ell +1 -n$.
On the other hand the scheme $2P$ is well poised and the  $\CC$-vector space generated
by the $E_i$ for
$1\leqslant i\leqslant \ell+1-n$ has codimension $n-1 < n$ in $L\otimes_\RR\CC$. 
So there exists at least one embedding $\tau$ such that the $(\ell +1 -n)\times r$ matrix
\[\left((\partial E_i /  \partial x_j)(P_\tau)\right)_{1\leqslant i\leqslant \ell+1-n, \,  1\leqslant j\leqslant r}\]
has maximal rank $r$. There exist $r$ integers
$1\leqslant i_1 < i_2 < \dots < i_r \leqslant \ell +1-n$
such that the value at $P_\tau$ of the minor determinant 
\begin{equation}\label{eq:Phinf}
  \Phi  = \det \left(  (\partial E_{i_k} /  \partial x_j)\right)_{1\leqslant k,\,  j \leqslant r}\end{equation}
is non-zero for some $\tau$ and thus for all $\tau$ by Galois action.  

\begin{proposition}[Number fields have small incomplete intersection
    models]\label{prop:sm}
  Let $\KK$ be a number field of degree $n$ and root
  discriminant $\delta_\KK$ over $\QQ$.
  Let $d\geqslant 5$ and $r\geqslant 1$ be rational  integers such that
  \[n(r+1)\leqslant {d+r\choose d}.\]
  Let 
  \begin{equation}\label{eq:H}
    H(\KK,d,r) = \ell^{\ell/2n}
    \times {d+r\choose d}^{1/2}\left(n^2d(r+1)\delta_\KK^{2}\right)^{d}\end{equation} where
\[\ell = {d+r\choose d}-n.\]  
Then $\KK/\QQ$ has an incomplete intersection model in dimension
$r$, degree $d$, and height $H(\KK,d,r)$.
\end{proposition}

In order to prove Theorem \ref{th:smnf}
we specialize  the values of the parameters $r$ and $d$ in Proposition \ref{prop:sm}. We will take $d=r$. It is evident
that ${2r \choose r}\geqslant 2^r$ so
\[\frac{1}{r+1}{2r \choose r}\geqslant 2^{\frac{r}{2}}\]
for $r$ large enough.
Further  \[\frac{1}{r+2}{2r+2 \choose r+1}\leqslant \frac{1}{r+1}{2r \choose r}\times 4.\]
We choose $r$ to be the smallest positive integer such that
$n(r+1)\leqslant {2r\choose r}$. We have
\begin{equation}\label{eq:choo}
  n(r+1)\leqslant {2r\choose r}\leqslant 4n(r+1) \text{ and } r\leqslant 3\log n \end{equation}
for $n$ large enough.
We deduce that $\ell = {2r\choose r}-n\leqslant 4n(r+1)\leqslant \cQ n\log n$.
So \[\ell^{\ell/2n}\leqslant n^{\cQ \log n}.\]
From Equation (\ref{eq:choo}) we deduce
that ${2r\choose r}\leqslant \cQ n\log n$.
Also $n^2d(r+1)\leqslant \cQ n^2\log^2 n$ and
\begin{equation}\label{eq:HK2}
  \left(n^2d(r+1)\right)^{d}\leqslant n^{\cQ \log n}.\end{equation}
  So   the coefficients of equations
$E_i$ are bounded in absolute value by
\begin{equation}\label{eq:HK}
  H(\KK,r,r)\leqslant (n\delta_\KK)^{\cQ \log n}.\end{equation}
This proves  Theorem \ref{th:smnf}.

\subsection{Interpolants for number fields}\label{sec:interpolnf}

In order to
prove Theorem \ref{th:sintnf}
we first look for  a prime integer $p$
such that the value at $P_\tau$ of the determinant $\Phi$ of Equation (\ref{eq:Phinf})
is prime to  $p$.
Let $\Psi  \in \ZZ$
be the product of all $\Phi (P_\tau)$
for all embeddings $\tau$. We first bound $\Psi$
in absolute value.
For every ${1\leqslant k,\,  j \leqslant r}$
the polynomial \[\partial E_{i_k} /  \partial x_j
\in \ZZ[x_1, \ldots, x_r]\] has  total degree
$\leqslant  r$ and height $\leqslant  rH(\KK,r,r)$.
Using Lemma \ref{lemma:eval},
Equation (\ref{eq:HK})
and Equation (\ref{eq:HK2})
we deduce
that the evaluation of
$\partial E_{i_k}/\partial x_j$ at $P_\tau$
is bounded in absolute value by
\begin{equation}\label{eq:G}
  G=\left(n^2r(r+1)\delta_\KK^2\right)^rrH(\KK,r,r){2r\choose r}\leqslant (n\delta_\KK)^{\cQ\log n}\end{equation}
We deduce that
\[|\Phi (P_\tau)| \leqslant G^rr!\leqslant (n\delta_\KK)^{\cQ(\log n)^2}
\text{\,\, and \,\,}
|\Psi | \leqslant (n^n|d_\KK|)^{\cQ(\log n)^2}.\]
Using an estimate \cite[Chapter I \S 2.6, Corollary 10.1]{Tenenbaum}
for the Tchebychev function
$\theta$  we deduce that there exists a prime \[p\leqslant
\cQ  n(\log n+\log \delta_\KK)(\log n)^2\] that does not divide
$\Psi$.

The ring $\ZZ[\kappa_1, \ldots, \kappa_r]$
is an order in $\KK$ which is unramified
at $p$. Let $\pp$ be  a place above
$p$. Let $s$ be 
the degree of inertia at $\pp$.
We set $q=p^s$ and we
fix an isomorphism between $\ZZ_q$ and the completion
of $\ZZ[\kappa_1, \ldots, \kappa_r]$  at $\pp$.
We call
$b_1^\infty$, \ldots, $b_r^\infty$ the images
of $\kappa_1$, \ldots, $\kappa_r$ by this bijection.
We let $m$ be the smallest positive integer such that
$q^m=p^{ms}$  is bigger than  $G^n$ where $G$ is given in Equation (\ref{eq:G}).
We check  that \[ms\log p\leqslant n\log G+s\log p \leqslant n(\log G+\log p)
\leqslant \cQ n\log n (\log n+\log\delta_\KK).\]

We set $b_j^m = b_j^\infty \bmod p^m \in  \ZZ_q/p^m$ for every
$1\leqslant j\leqslant r$.
The point \[b^m=(b_1^m, \ldots, b_r^m)\] has an interpolant in degree $d=r$ and height $H(\KK,r,r)$ and this interpolant is
$b^\infty$.
Indeed let $E(x_1, \ldots, x_r)$ be a polynomial
vanishing at $b^m$ and
with total degree $\leqslant r$
and height $\leqslant H(\KK,r,r)$.
The evaluation
of $E$
at $P_\tau$ is bounded from above
by $G$ according to Lemma \ref{lemma:eval}.
Its norm is thus an integer bounded in absolute value
by $G^n$ and divisible by $p^{ms}>G^n$. So $E$ vanishes at $P_\tau$
and at $b^\infty$ just as well.
Unicity follows from the existence of $r$ equations
with degree $\leqslant r$ and
heigth $\leqslant H(\KK,r,r)$ having a non-zero Jacobian
determinant modulo  $\pp$.
This  finishes the proof of Theorem \ref{th:sintnf}.

\section{Conclusion}

We have constructed short models of two kinds
for a global   field: as  an irreducible
component of some complete intersection or as an interpolant
of a finite affine scheme. 
Theorems \ref{th:siiff} and
\ref{th:sintff} on the one hand and Theorems
\ref{th:smnf} and \ref{th:sintnf} on the other
hand present evident
similitudes. The awaited factor $g+n(1+\log_qn)$
in the geometric situation has counterpart
$n(\log n+\log \delta_\KK)$ in the arithmetic
situation. The not so welcome factor
$(\log  n)^3$ or
$(\log  n)^2$ is the same in either cases.
The methods are parallel to some extent.
Interpolation plays a crucial role on either
sides, especially the theorem of
Alexander
and   Hirschowitz. We do not use the geometric
equivalent of Minkowski space however although it
exists, see \cite{rosen}. The reason is that
basic arguments about vector bundles over
$\PP^1$ suffice for our purpose. 
Another difference is that in the arithmetic
case our constructions are consistent with the
best known bounds for the number of number
fields of given degree and bounded discriminant
\cite{jmc2019,LT,lee}.  In the geometric
situation  however, cohomological methods produce
much better estimates \cite{JK}, but no effective construction.

\bibliographystyle{plain}
\bibliography{short}
\end{document}